\documentclass{birkart}

 \newtheorem{thm}{Theorem}[section]

 \newtheorem{prop}[thm]{Proposition}
 \theoremstyle{definition}
 \newtheorem{defn}[thm]{Definition}
 \theoremstyle{remark}
 \newtheorem{rem}[thm]{Remark}
 
 \numberwithin{equation}{section}

\begin{document}
%
\title[Polynomial Approximation]
 {Regular polynomial interpolation and approximation of global solutions of linear partial differential equations}
\author{ J\"org Kampen}
\address{%
Weierstrass Institute
for Applied Analsis and
Stochastics
Mohrenstrasse 39\\
10117 Berlin\\
Germany}
\email{kampen@wias-berlin.de}

\thanks{This work was completed with the support of DFG (Matheon)}



\date{Mai 31, 2007}

\begin{abstract}
We consider regular polynomial interpolation algorithms on recursively defined sets of interpolation points which approximate global solutions of arbitrary well-posed systems of linear partial differential equations.
 Convergence of the 'limit' of the recursively constructed family of polynomials to the solution and error estimates are obtained from a priori estimates for some standard classes of linear partial differential equations, i.e. elliptic and hyperbolic equations.  
Another variation of the algorithm allows to construct polynomial interpolations which preserve systems of linear partial differential equations at the interpolation points. We show how this can be applied in order to compute higher order terms of WKB-approximations of fundamental solutions of a large class of linear parabolic equations. The error estimates are sensitive to the regularity of the solution. Our method is compatible with recent developments for solution of higher dimensional partial differential equations, i.e. (adaptive) sparse grids, and weighted Monte-Carlo, and has obvious applications to mathematical finance and physics.
\end{abstract}

\maketitle
\section{Introduction}
This work shows how multivariate interpolation techniques can be combined with analytic information of linear partial differential equations (i.e.  a priori estimates and/or WKB representations of solutions) in order to design efficient and accurate numerical schemes for solving (systems) of linear partial differential equations. These schemes are nothing but sequences of multivariate polynomials which are constructed recursively such that they solve a given linear system of partial differential equations on a finite discrete set of interpolation points. However, additional information is needed in order to ensure that the sequence of interpolation polynomials  converges to a (or, if uniqueness is proved, the) global solution of a given linear system of partial differential equations. As we shall see, this information can be provided by a priori estimates which in turn lead us to error estimates in regular norms dependent on the regularity of the solution. We examine the situation in the case of linear elliptic equations with variable coefficients. Another possibility is that (more or less) explicit representations of solutions are known which lead to problems which are easier to solve. A prominent example is the WKB-expansion which was investigated in \cite{Ka}. The recursive structure of WKB coefficient functions and the error analysis lead us to the problem of regular polynomial approximation. In this introductionary Section we our method on an abstract level.
 
\subsection{Regular polynomial interpolation}

Since we are interested in the relationship between multivariate polynomial interpolation and approximation of solutions of partial differential equations, our focus will be on multivariate polynomial interpolation. However, in order to make basic ideas more accessible we shall describe algorithms in the univariate case first and then generalize to the multivariate case. It is well known that polynomial interpolation in the multivariate case is quite different from the univariate case in general. However, in our approach which aims at solving linear systems of partial differential equations or aims at supplementing certain strategies of solving partial differential equations many features are already present in the univariate framework. In order to avoid misunderstandings, we dwell a little on this point.
Classically, the problem of multivariate interpolation can be stated as follows (cf. \cite{Sau2}):

{\it
Given a set of interpolation points $\Theta=\left\lbrace x_1,\cdots ,x_N \right\rbrace $  and an N-dimensional

 space $P_{\Theta}$ of polynomials  find, for given values $y_1,\cdots,y_N$, a unique polynomial $f\in P$ such that
\begin{equation}
f(x_j)=y_j,~j\in 1, \cdots ,N.  
\end{equation}\it}
In this form it turns out that there is an intricate relation between sets of interpolation points and interpolation spaces that must be satisfied in order that the problem can be considered to be well-posed. Either we have to make some restrictions concerning the set of interpolation points $\Theta$ (cf. \cite{Sau2}) or we consider $\Theta$ to be fixed and consider the problem of constructing the polynomial space $P_{\theta}$ (cf.\cite{boor}). This amounts to a construction of the map 
\begin{equation}
\Theta \rightarrow P_{\Theta}
\end{equation}
with additional constraints such as minimality of degree (cf. \cite{Sau2, boor}) or monotonicity (cf. \cite{boor}). In this paper we are interested in interpolation algorithms with the following features
\begin{itemize}

\item there are no essential restriction on the discrete set $\Theta$ of interpolation points except that $\Theta \subset D$, where $D$ is the domain of the function to be interpolated.  

\item the map $\Theta \rightarrow P_{\Theta}$ is monoton (indeed our basic algorithm is an extension of multivariate versions of Newton's interpolation algorithm).

\item the algorithm can be extended to vector valued interpolation functions $g: D\subseteq {\mathbb R}^n\rightarrow {\mathbb R}^k$ and if $g$ satisfies a system of linear partial differential equations, then the interpolation polynomial $p$ solves the same system of linear partial differential equations on the given set $\Theta$ of interpolation points.   

\item the algorithm is numerically stable and practical with respect to the problem that the interpolation function $f$ and arbitrary set of partial derivatives of $f$ are to be interpolated simultaneously. For the application of higher order approximation of the fundamental solution of linear parabolic equations we comute accurate approximations of derivatives of smooth functions up to order $10$ in order to obtain an approxmation of order $5$ of the WKB-expansion of the fundamental solution.

\item the algorithm can be refined in order to solve well-posed linear systems of partial differential equations directly.

\item the algorithm can be combined with collocation methods in an efficient way; it can be partially parallelized.   

\item the algorithm allows for error estimates which depend on the regularity of the solution such that the algorithm is compatible with methods for higher dimensional problems of linear systems of partial differential equations such as sparse grids, adaptive sparse grids, and weighted Monte-Carlo.
\end{itemize}

First we consider the problem of polynomial approximation $p$ of a regular (i.e smooth or finitely many times differentiable) function  
\begin{equation}
f:D\subseteq {\mathbb R}^n\rightarrow {\mathbb R}
\end{equation}
defined on discrete subset of $\Theta\subset D$ where for $m$ given linear partial differential operators
\begin{equation} 
L_i=\sum_{|\alpha| \leq \beta_i} a^i_{\alpha}(x)\partial_{\alpha},
\end{equation}
we require that
\begin{equation}
L_i f(x_j)=L_i p(x_j) \mbox{ for } 1\leq i \leq m
\end{equation}
for some finite set of points $x_j\in \Theta\subset D$. As indicated above we shall allow that the interpolation set $\Theta$ can be constructed recursively (and, hence, extended arbitrarily within the domain of the interpolation function).
Investigations of specific instances of this problem can be found in the literature  on polynomial interpolation (cf. the survey paper of \cite{Sau} for the development up to the year 2001). Note that other algorithms of natural interpolation of $C^k$-functions have been proposed (cf.\cite{Gasca} for hints at the history and further references).

The paper is organized as follows. In Section 1.2 we introduce the partial differential equations for which we seek global regular interpolation polynomials of their global solutions. All basic types of partial differential equations, i.e. elliptic equations, parabolic equations, and hyperbolic equations are considered. While the basic algorithm is quite similar for each type of partial differential equation, we shall see, however, that the convergence of the scheme of recursively defined interpolation polynomials depends on very different a priori estimates for different type of equations. In case of second order elliptic equations classical Schauder boundary estimates can be used, while in the case of hyperbolic equations energy estimates are considered. In the case of parabolic equations we refer back to Safanov-Krylov estimates considered in the context of the truncation error analysis of WKB-expansions. In Section 2.1 we introduce first an extension of Newton's polynomial algorithm which interpolates a given function and its derivatives up to some given order $k$ simultaneously. Section 2.2. describes  a variation of this algorithm which interpolates a given function such that a given set of partial differential equations is preserved. Section 3 discusses the extension to the multvariate case. In Section 4 we refine the algorithm and construct polynomials which satisfy a given linear (i.g. partial) differential equation on a given set of interpolation points, i.e. there is no given function to be interpolated. In Section 5 we consider refinements which show how polynomials constructed on disjoint sets of interpolation points can be synthesized in order to get one polynomial which interpolates on the union of sets of interpolation points. Naturally, parallelization is consideredin this context. In Section 6 we show how a priori estimates of elliptic equations (standard Schauder boundary estimates) and hyperbolic equations (energy estimates) lead to convergent schemes implied by error estimates. Section 7 discusses a special use of regular polynomial interpolation for parabolic equations where the global solution is given in the form of a WKB- expansion. Section 8 provides a numerical example of global regular polynomial interpolation of a locally analytic function up to the third derivative. In Section 9 we provide a summary and give an outlook on current research and research in the near future. Before we start with the description of the algorithm, we state the typical linear partial differential equations and indicate the different types of approximations and error estimates which we aim at.

\subsection{Regular interpolation and partial differential equations}
 We consider the three standard types of linear partial differential equations, namely elliptic equations, parabolic equations, and hyperbolic equations, and exemplify different types of application and extension.

\begin{itemize}

\item The most popular examples of elliptic partial differential equations are of the second order form, i.e.

\begin{equation}
\sum_{j, k}^na_{jk}(x)\frac{\partial^{2}u}{\partial x_j \partial x_k}+\sum_{l}b_l(x)\frac{\partial u}{\partial x_l }+c(x)u=f(x),
\end{equation}
to be solved on a domain $\Omega\subseteq {\mathbb R}^n$ with the
boundary condition
\begin{equation}
u\mbox{\Big |}_{\partial \Omega}=g
\end{equation}
for some function 
$f:{\partial \Omega}\rightarrow {\mathbb R}$ which is usually assumed to be Lipschitz continuous at least.
Here,  $a_{jk}$ are (at least) measurable coefficient functions satisfying for some constant $c$, and ellipticity means that
\begin{equation}
\sum_{jk}a_{jk}(x)\xi_i\xi_j\geq c>0 \mbox{ (uniformly in $x$)}.
\end{equation}
We construct an extension of the polynomial interpolation algorithm which produces a multivariate polynomial solving this elliptic equation on an arbitrary grid of interpolation points. In order to obtain error estimates b standard boundary Schauder estimates in this paper we shall make some regularity assumptions. We derive convergence of the family of multivariate polynomials constructed by our our interpolation scheme to the global solution of the linear elliptic equation on a bounded domain and we derive error estimates from a priori estimates.

\item Parabolic equations of the form
\begin{equation}\label{parab}
\frac{\partial u}{\partial  t}-Lu=0,
\end{equation}
on $D:=\Omega\times (0,T)$, ($\Omega\subseteq {\mathbb R}^n$, with
\begin{equation}
u(0,x)=\delta_y(x):=\delta(x-y),~y\in{\mathbb R}^n,
\end{equation}
where $\delta$ is the Dirac delta distribution, and where
\begin{equation}
Lu\equiv \frac{1}{2}\sum_{ij}a_{ij}(x)\frac{\partial^2 u}{\partial x_i\partial x_j}+\sum_i b_i(x)\frac{\partial u}{\partial x_i}
\end{equation}
is an elliptic operator. The solution of this equation is called fundamental solution, because solutions of standard parabolic initial-value boundary problems can be represented by convolution integrals of data functions with the fundamental solution. The standard assumptions for such a fundamental solution to exist are
\begin{itemize}
\item[(A)] The operator $L$ is uniformly parabolic in ${\mathbb R}^n$, i.e. there exists $0<\lambda <\Lambda <\infty$
such that for all $\xi\in {\mathbb R}^n\setminus \{0\}$
$$
0<\lambda |\xi|^2\leq \sum_{i,j=1}^n a_{ij}(x)\xi_i\xi_j\leq \Lambda |\xi|^2.
$$

\item[(B)] The coefficients of $L$ are bounded functions in ${\mathbb R}^n$ which are uniformly
H\"older continuous of exponent $\alpha$ ($\alpha \in (0,1)$).
\end{itemize}
If some regularity assumptions on the coefficients hold in addition, then it can be shown that the fundamental solution $p$ is of the form
\begin{equation}\label{frep}
 p(t,x,y)=\frac{1}{\sqrt{2\pi  t}^n}\exp\left( -\frac{d^2(x,y)}{2 t}+\sum_{k\geq 0}c_k(x,y)t^k\right),
 \end{equation}
with some regular coefficient functions $d^2$ and $c_k$. We shall show how our regular polynomial interpolation algorithm can be used to compute the fundamental solution in terms of this representation. 

\begin{rem}
The algorithm designed in the case of elliptic equations can  be applied to the parabolic case directly, of course. However, it turns out that the convergence is better if the special representation \eqref{frep} is used. 
\end{rem}

\item 
As an example of a hyperbolic equation we consider an equation of the form
\begin{equation}
Lu=f \mbox{ in } \Omega, 
\end{equation}
where
\begin{equation}
Lu\equiv \sum_{ij} h_{ij}\frac{\partial u}{\partial x_i\partial x_j}+\sum_i\frac{\partial}{\partial x_j}+c(x)u
\end{equation}
and $(h_{ij})$ is a symmetric matrix of signature $(n,1)$, if $\mbox{dim}\Omega =n+1$. We assume that some $O\subset \Omega$ is bounded by two spacelike surfaces $\Sigma_i$ and $\Sigma_e$ and swept out by a family of spacelike surfaces $\Sigma_e(s)$. We assume the initial conditions
\begin{equation}
u=g \mbox{ and } du=\omega 
\end{equation}
where $g$ is a function on $\Omega$ and $\Omega$ is a 1-form.

\end{itemize}

\section{Interpolation algorithm (univariate case)}

We start with the description of the algorithm which produces polynomials which satisfy some given requirements on interpolation points. Our starting point is an extension of Newton's polynomial interpolation method such that the interpolation polynomial and its derivatives up to a given order $k$ (an integer) equal a given function and its derivatives up to order $k$ at the interpolation points. For simplicity of representation and since the essential features of the algorithm can be demonstrated for one dimensional functions, we describe our ideas first in the univariate case and then generalize to the multivariate case in the next section. 

\subsection{Extension of Newton's method}
Let us recall the Newtonian interpolation for an univariate function
\begin{equation}
\begin{array}{ll}
f:[a,b]\subset {\mathbb R}\rightarrow {\mathbb R}.
\end{array}
\end{equation}
Given a discrete set of interpolation points $D=\{x_0,x_1\cdots ,x_N\}\subset [a,b]$ we want to construct a polynomial 
\begin{equation}
\begin{array}{ll}
p:[a,b]\subset {\mathbb R}\rightarrow {\mathbb R} \mbox{ such that}\\
\\
f(x_i)=p(x_i) \mbox{ for all } x_i \in D.
\end{array}
\end{equation}
The idea of the basic Newton interpolation algorithm is that instead of looking for some polynomial of form $\sum_{i=1}^N b_i x^i$ for some constants $b_i$ we may write
\begin{equation}
\sum_{l=0}^N a_l \Phi_l(x)
\end{equation}
with 
\begin{equation}
\Phi_0(x)=1 \mbox{ and } \Phi_l(x)=
\Pi_{i=0}^{l}(x-x_i) \mbox{ for $l\geq 1$.}
 \end{equation}
In order to determine $a_0,\cdots a_{N}$ we then may solve the system
\begin{equation}
R_0a:=\left[ \begin{array}{ccccccc}
&1 &0 & 0 &\cdots &0\\
&1 &\phi_1(x_1) & 0&\cdots &0\\
&1 & \phi_1(x_2) &\phi_2(x_2) &\cdots & 0\\
&\vdots & \vdots &\vdots &\vdots\\
&1 & \phi_1(x_{N}) &\phi_2(x_{N}) &\cdots & \phi_N(x_{N})
\end{array}
  \right] \left[ \begin{array}{ccccccc}
a_0\\
a_1\\
a_2\\
\vdots\\
a_{N} 
\end{array}
 \right] =\left[ \begin{array}{ccccccc}
f(x_0)\\
f(x_1)\\
f(x_2)\\
\vdots\\
f(x_{N})
\end{array}
\right]
\end{equation}
This leads to an $L^2$-approximation of the function $f$ similar to the Gaussian algorithm. Note however, that the matrix $R_0$  is a lower diagonal. Hence the linear system can be solved easily. Moreover the matrix condition number is much better than that of the Vandermonde matrix used in the classical Gaussian interpolation. We extend this idea to a $C^k$-norm interpolation, i.e. we design an algorithm that approximates $f$ up to the $k$-th derivative, i.e. we construct a polynomial 
\begin{equation}
\begin{array}{ll}
q:[a,b]\subset {\mathbb R}\rightarrow {\mathbb R} \mbox{ such that}\\
\\
f^{(l)}(x_i)=q^{(l)}(x_i) \mbox{ for all } x_i \in D \mbox{ and all } l\leq k,
\end{array}
\end{equation}
where for a function $g:[a,b]\subset {\mathbb R}\rightarrow {\mathbb R}$ $g^{(l)}$ denotes the derivative of order $l$ while $g=g^{0}$. We consider the polynomial
\begin{equation}
\sum_{m=0}^{(N+1)(k+1)-1}a_m\Phi_{m,k}(x) 
\end{equation}
where
\begin{equation}
\Phi_{m,k}(x)=(x-x_{m \mbox{ div}(k+1)})^{m \mbox{ mod}(k+1)}\Pi_{l=0}^{m \mbox{div}(k+1)-1}(x-x_l)^{k+1},
\end{equation}
where, by convention, we understand
\begin{equation}
\Pi_{l=0}^{-1}(x-x_l)^{k+1}:=1.
\end{equation}
For simplicity of notation we sometimes use the abbreviations
\begin{equation}
p(m)= m  \mbox{div} (k+1) \mbox{ and } q(m)=m\mbox{mod}(k+1).
\end{equation}
Next we define
\begin{equation}
\Phi^{(l)}_{m,k}(x):= \frac{d}{dx^l}\Phi_{m,k}(x),
\end{equation}
and for each $k\geq 1$ the linear system
\begin{equation}
R_k\left[ \begin{array}{ccccccc}
a_0\\
a_1\\
a_2\\
\vdots\\
a_{(k+1)(N+1)-1} 
\end{array}
 \right] =\left[ \begin{array}{ccccccc}
f(x_0)\\
f'(x_0)\\
\vdots\\
f^{(k)}(x_0)\\
f(x_1)\\
\vdots\\
f^{(k)}(x_{(k+1)(N+1)-1})
\end{array}
\right]
\end{equation}
where $R_k$ is a $(N+1)(k+1)\times (N+1)(k+1)$-matrix determined by $(k+1)\times (k+1)$ matrices $A_k^{lm}$ as follows:
\begin{equation}
R_k:=\left[ \begin{array}{ccccccc}
&A^{00}_k &Z_k & Z_k &Z_k &\cdots &Z_k\\
&A^{10}_k &A^{11}_k  & Z_k & Z_k &\cdots &Z_k\\ 
&A^{20}_k &A^{21}_k & A^{31}_k & Z_k &\cdots &Z_k\\
&\vdots &\vdots & \vdots &\vdots &\vdots &\vdots\\  
&A^{N0}_k & A^{N1}_k & A^{N2}_k &A^{N3}_k &\cdots &A^{NN}_k\\ 
\end{array}
\right],
\end{equation}
where $Z_k$ is the $(k+1)\times (k+1)$ matrix with $0$ entries, and  
\begin{equation}
A^{ij}_k=A^i_k(x_j)
\end{equation}
with
\begin{equation}
A^{ij}_k:=\left[ \begin{array}{ccccccc}
&\Phi_{(k+1)p(i),k}(x_j) &\Phi_{(k+1)p(i)+1,k}(x_j) & \Phi_{(k+1)p(i)+2,k}(x_j) &\cdots &\Phi_{(k+1)p(i)+k,k}(x_j)\\
&\Phi^{(1)}_{(k+1)p(i),k}(x_j) &\Phi^{(1)}_{(k+1)p(i)+1,k}(x_j) & \Phi^{(1)}_{(k+1)p(i)+2,k}(x_j) &\cdots &\Phi^{(1)}_{(k+1)p(i)+k,k}(x_j)\\

&\Phi^{(2)}_{(k+1)p(i),k}(x_j) &\Phi^{(2)}_{(k+1)p(i)+1,k}(x_j) & \Phi^{(2)}_{(k+1)p(i)+2,k}(x_j)  &\cdots &\Phi^{(2)}_{(k+1)p(i)+k,k}(x_j)\\

&\vdots &\vdots & \vdots &\vdots  &\vdots\\  
&\Phi^{(k)}_{(k+1)p(i),k}(x_j) &\Phi^{(k)}_{(k+1)p(i)+1,k}(x_j) & \Phi^{(k)}_{(k+1)p(i)+2,k}(x_j) &\cdots &\Phi^{(k)}_{(k+1)p(i)+k,k}(x_j)
\end{array}
  \right].
\end{equation}
Note that
\begin{equation}
A^{00}_k:=\left[ \begin{array}{ccccccc}
&1 &0 & 0 &0 &\cdots &0\\
&0 &1 & 0 &0 &\cdots &0\\ 
&0 &0 & 2 &0 &\cdots &0\\
&\vdots &\vdots & \vdots &\vdots &\vdots &\vdots\\  
&0 &0 & 0 &0 &\cdots &k!
\end{array}
  \right]. 
\end{equation}
This leads to a system which can be solved row by row. It is therefore very easy to implement and numerically well-conditioned.

\begin{rem}
In order to avoid large entries in the matrices $A^{lm}_k$ one may consider basis functions of form $\frac{1}{l!}\Phi^{(l)}_{(k+1)p(i),k}$, but we do not deal with the peculiar niceties of computation here.
\end{rem}
\subsection{Interpolation preserving linear systems of differential equations}
The preceding algorithm can be adapted it in order to construct a polynomial approximation $p$ of $f$ where the $k$ differential operators
\begin{equation}
L_if(x)=\sum_{j \leq q_i} a^i_{j}(x)\frac{d}{d x^j}f(x), i=1,\cdots,k
\end{equation}
are preserved on a discrete set of points $\Theta=\{x_0,\cdots ,x_N\}$ in the sense that
\begin{equation}
L_if(x_j)=L_ip(x_j) \mbox{ for } x_j\in \Theta.
\end{equation}
At this point the linear system of the operators $\left\lbrace L_i|1\leq i\leq k\right\rbrace $ is quite arbitrary; we just assume that the operators are defined pointwise, i.e. $x\rightarrow a^i_j(x)$ are classical functions which can be evaluated pointwise (at least on the set of interpolation points). Note that we do not ask about convergence of a family of interpolation polynomials to at this point. There are several possibilities to extend our preceding algorithm. One is the following. 
Let 
\begin{equation}
Q_i:=\left\lbrace j|a^i_j\neq 0\right\rbrace 
\end{equation}
and define 
\begin{equation}
L_i^m=\sum_{j\in Q_i, j\leq m}a^1_{i_j}(x)\frac{d^{i_j}}{d x^{i_j}}.
\end{equation}
We start with 
\begin{equation}
Q_1=\left\lbrace i_{11},\cdots ,i_{1r_1}\right\rbrace ,
\end{equation}
and assume that
\begin{equation}
i_{11}< \cdots < i_{1r_1}
\end{equation}
We consider first the interpolation point $x_0$ and start with the following ansatz for the interpolation polynomial
\begin{equation}
p_{10}(x)=\sum_{i_{1j}\in Q_1}b^{10}_{i_j}(x-x_{0})^{i_{1j}}. 
\end{equation}
We assume $f(x_0)=p_{10}(x_0)=0$ w.l.o.g. ; we shall see later how we interpolate values of $f$ different from zero at the other interpolation points $x_1,\cdots, x_N$. First we apply the operator
\begin{equation}
L_1^{i_1}\equiv a^1_{i_1}(x)\frac{d^{i_1}}{d x^{i_1}}
\end{equation}
to $f$ and $p_{10}$ at $x_0$. This leads to 
\begin{equation}
i_1!b^{10}_{i_1}=a^1_{i_1}(x_0)\frac{d^{i_1} f}{d x^{i_1}}(x_0)~~\Rightarrow~~b^{10}_{i_1}=\frac{1}{i_1!}a^1_{i_1}\frac{d^{i_j} f}{d x^{i_1}}(x_0 )
\end{equation}
Inductively we assume that the coefficients $b^{10}_{i_j}$ have been defined up to the index $i_m$ for some $m< r_1$ and that the operator $L_1^{i_m}$ has been defined accordingly. We apply the operator 
\begin{equation}
L^{i_{m+1}}_1\equiv L^{i_{m}}_1 + a^1_{i_{m+1}}(x)\frac{d^{i_{m+1}}}{d x^{i_{m+1}}}
\end{equation}
to $f$ and $p_{10}$ at $x_0$. For an integer $s$ with $m+1\leq s\leq r_1$ define
\begin{equation}
p_{10}^s(x)=\sum_{j=1}^{s}b^{10}_{i_j}(x-x_0)^{i_j}.
\end{equation}
Then we have
\begin{equation}
\begin{array}{ll}
L^{i_{m+1}}_1p_{10}(x_0)= L^{i_m}_1p_{10}(x_0) +a^1_{i_{m+1}}(x_0)\frac{d^{i_{m+1}}}{d x^{i_{m+1}}}p_{10}(x_0)=\\
\\
L^{i_m}_1p_{10}^{i_m}(x_0)+i_{m+1}! a^1_{i_{m+1}}(x_0)b^{10}_{i_{m+1}}=L^{i_{m+1}}_1f(x_0).
\end{array}
\end{equation}
This gives $b^{10}_{i_{m+1}}$.
Next inductively assume that an interpolation polynomial $p_{1k}$ has been constructed which interpolates $f$ on the set of interpolation points $\left\lbrace x_0,\cdots ,x_k\right\rbrace $ for some positive integer $k$ with $k< N $ subject to the condition
\begin{equation}
L_1 f(x_i)=p_{1k}(x_i) \mbox{ for } 1\leq i\leq k. 
\end{equation}
First we extend that polynomial in order to interpolate $f$ at the point $x_{k+1}$. We consider the ansatz
\begin{equation}
p^0_{1(k+1)}(x)=p_{1k}(x)+b^{1(k+1)}_0\Pi_{l=0}^k(x-x_l)^{q_1}.
\end{equation}
We then get $b^{1(k+1)}_0$ from the equation
\begin{equation}
p^0_{1(k+1)}(x_{k+1})=f(x_{k+1}).
\end{equation}
The ansatz for $p_{1(k+1)}$ (i.e. the interpolation polynomial which preserves $L_1 f$ on the set of interpolation points $\left\lbrace x_1,\cdots ,x_{k+1}\right\rbrace$) is
\begin{equation}
p_{1(k+1)}(x)=p^0_{1(k+1)}(x)+\sum_{i_j\in Q_1} b^{1(k+1)}_{i_j}(x-x_{k+1})^{i_j}\Pi_{l=0}^k(x-x_l)^{q_1+1}
\end{equation}
and the determination of coefficient constants $b^{1(k+1)}_{i_j}$ is similar to the procedure for the interpolation point $x_0$ described above.
Proceeding inductively, we are lead to the polynomial $p_{1}$ which interpolates $f$ at the interpolation points of $\Theta =\left\lbrace x_0,\cdots ,x_N\right\rbrace $ such that
\begin{equation}
L_1 p_1(x_j)= L_1f(x_j) \mbox{ for all } x_j\in \Theta.
\end{equation}
Finally assuming that for some integer $s<k$ the polynomial $p_s$ satisfies the condition that
\begin{equation}
\begin{array}{ll}
p_s(x_j)=f(x_j) \mbox{ for } x_j\in \Theta\\
\\
L_ip_s(x_j)=L_if(x_j) \mbox{ for } x_j\in \Theta \mbox{ and } i\leq s,
\end{array}
\end{equation}
it is clear that we only need to consider the reduced operator
\begin{equation}
L_{s+1}\equiv \sum_{i_jj\in Q_{s+1}\setminus \cup_{i=1}^{s} Q_i}a^{s+1}_{i_j}(x)\frac{d^{i_j}}{d x^{i_j}}.
\end{equation}
and proceed analogously.

\section{Extension to the multivariate case}

Next we consider generalizations to the multivariate case. There are several possibilities but the most simple seems to be the following. First we formulate the problem in a way that will turn out to be useful in the context of polynomial interpolation of global solutions of linear systems of partial differential equations. In its most simple form it is a form of multivariate Newton interpolation: given a function 
\begin{equation}
\begin{array}{ll}
f: S\subset {\mathbb R}^n\rightarrow {\mathbb R} 
\end{array}
\end{equation}
we want to construct a polynomial
\begin{equation}
\begin{array}{ll}
p:S\subset {\mathbb R}^n\rightarrow {\mathbb R} \mbox{ such that}\\
\\
f(x_i)=p(x_i) \mbox{ for all } x_i \in D \subseteq S,
\end{array}
\end{equation}
where $D=\left\lbrace x_0,x_1,\cdots ,x_n\right\rbrace $ is some discrete sets of points in ${\mathbb R}^n$ whose coordinates will be denoted by superscript indices as $x_i^j,~j=1,\cdots ,n$. This is done then by recursive definition of polynomials $p_0, p_1,\cdots $. First, define 
\begin{equation}
\begin{array}{ll}
p_0(x)\equiv f(x_0).
\end{array}
\end{equation}
Next, ansatz and equation
\begin{equation}
\begin{array}{ll}
p_1(x)\equiv f(x_0)+a_1\Pi_{i=1}^n(x^i-x_0^i)=f(x_1)
\end{array}
\end{equation}
leads to the determination of $p_1$ by 
\begin{equation}
\begin{array}{ll}
a_1=\frac{f(x_1)-f(x_0)}{\Pi_{i=1}^n(x^i-x_0^i)}
\end{array}
\end{equation}
Next assume that $p_0,p_1,\cdots, p_q$ have been defined. Then ansatz and equation
\begin{equation}
\begin{array}{ll}
p_{q+1}(x_{q+1})\equiv p(x_{q+1})+a_{q+1}\Pi_{k=0}^{q}\Pi_{i=1}^n(x^i-x_k^i)=f(x_{q+1})
\end{array}
\end{equation}
leads to the determination of $p_{q+1}$ by
\begin{equation}
\begin{array}{ll}
a_{q+1}=\frac{f(x_{q+1})-p_q(x_{q+1})}{\Pi_{k=0}^q\Pi_{i=1}^n(x^i-x_k^i)}
\end{array}
\end{equation}

%

\subsection{Extension of Newton's method}
Next we extend a multivariate version of Newton's method, i.e. we design an algorithm that approximates $f$ up to the $\beta$-th derivative ($\beta=(\beta_1,\cdots ,\beta_n)$ being some multiindex) where we construct a polynomial 
\begin{equation}
\begin{array}{ll}
q: S\subset {\mathbb R}\rightarrow {\mathbb R} \mbox{ such that}\\
\\
\frac{\partial f}{\partial x^{\gamma}}(x_i)
=\frac{\partial q}{\partial x^{\gamma}}(x_i) 
\mbox{ for all } 
x_i \in D\subseteq S \mbox{ and all } \gamma \leq \beta.
\end{array}
\end{equation}
where $\beta$ is given (i.e. the multivariate substitute for $k$ in the univariate case described above), and ordering is in the following sense:
 \begin{defn}
 Let $x^{\alpha}$ and $x^{\beta}$ be monomials in ${\mathbb R}\left[x_1,\cdots,x_n\right]$. We say that $x^{\alpha}>x^{\beta}$ (  lexicographical order) if $\sum_{i}\alpha^i > \sum_i \beta^i$ or  $\sum_{i}\alpha^i=\sum_i \beta^i$, and in the difference $\alpha -\beta\in {\mathbb Z}^n$ the left-most non zero entity is positive. 
\end{defn}
Now, let $\alpha_0,\alpha_1,\cdots,\alpha_m,\cdots $ an enumeration of multiindices with respect to this ordering. We define a sequence of polynomials $p_{\alpha_0}, p_{\alpha_1},\cdots ,p_{\alpha_m},\cdots $ recursively. First, let
\begin{equation}
p_{\alpha_0}(x)=a_{\alpha_0}+\sum_{\gamma \leq \beta}a_{\alpha_0\gamma}\Pi_{i=1}^n(x^i-x^i_{\alpha_0})^{\gamma_i}.
\end{equation}
If $p_{\alpha_0},\cdots ,p_{\alpha_{m-1}}$ have been defined, then we define
\begin{equation}
\begin{array}{ll}
p_{\alpha_m}(x)=p_{\alpha_{m-1}}(x)+\\
\\
\sum_{\gamma \leq \beta}a_{\alpha_{m-1}\gamma}\Pi_{i=1}^n(x^i-x^i_{\alpha_{m-1}})^{\gamma^i}\Pi_{j=0}^{m-1}\Pi_{i=1}^n(x^i-x^i_{\alpha_j})^{\beta^i+1}.
\end{array}
\end{equation}
This leads to a linear system to be solved for a vector $\left(a_{\alpha_0},\cdots , a_{\alpha_N\beta} \right)$ of length $(N+1)\left(\sum_i \beta^i+1\right)$ 

\begin{equation}
R_{\beta}\left[ \begin{array}{ccccccc}
a_{\alpha_0}\\
\vdots\\
a_{\alpha_0\beta}\\
a_{\alpha_1}\\
\vdots\\
a_{\alpha_N\beta} 
\end{array}
 \right] =\left[ \begin{array}{ccccccc}
f(x_{\alpha_0})\\
\vdots\\
f^{(\beta)}(x_{\alpha_0})\\
f(x_{\alpha_1})\\
\vdots\\
f^{(\beta)}(x_{\alpha_N})
\end{array}
\right]
\end{equation}
with
\begin{equation}
R_{\beta}:=\left[ \begin{array}{ccccccc}
&A^{00}_{\beta} &Z_{\beta} & Z_{\beta} &Z_{\beta} &\cdots &Z_{\beta}\\
&A^{10}_{\beta} &A^{11}_{\beta}  & Z_{\beta} & Z_{\beta} &\cdots &Z_{\beta}\\ 
&A^{20}_{\beta} &A^{21}_{\beta} & A^{31}_{\beta} & Z_{\beta} &\cdots &Z_{\beta}\\
&\vdots &\vdots & \vdots &\vdots &\vdots &\vdots\\  
&A^{N0}_{\beta} & A^{N1}_{\beta} & A^{N2}_{\beta} &A^{N3}_{\beta} &\cdots &A^{NN}_{\beta}\\ 
\end{array}
\right] 
\end{equation}
We abbreviate $\sum \beta = \sum_i (\beta^i+1)$ and defining $p(m)=m\div \sum \beta$
we have
\begin{equation}
A^{ij}_k:=\left[ \begin{array}{ccccccc}
&\Phi_{{\tiny \sum} \beta p(i),\beta}(x_j) &\Phi_{{\tiny \sum} \beta p(i)+\beta_1,\beta}(x_j) & \Phi_{{\tiny \sum} \beta p(i)+\beta_2,\beta}(x_j) &\cdots &\Phi_{{\tiny \sum} \beta p(i)+\beta,\beta}(x_j)\\
&\Phi^{(\beta_1)}_{{\tiny \sum} \beta p(i),\beta}(x_j) &\Phi^{(1)}_{{\tiny \sum} \beta p(i)+\beta_1,\beta}(x_j) & \Phi^{(\beta_1)}_{(k+1)p(i)+\beta_2,\beta}(x_j) &\cdots &\Phi^{(\beta_1)}_{{\tiny \sum} \beta p(i)+\beta,\beta}(x_j)\\

&\Phi^{(\beta_2)}_{{\tiny \sum} \beta p(i),\beta}(x_j) &\Phi^{(\beta_2)}_{{\tiny \sum} \beta p(i)+\beta_1,\beta}(x_j) & \Phi^{(\beta_2)}_{\sum \beta p(i)+\beta_2,k}(x_j)  &\cdots &\Phi^{(\beta_2)}_{\sum \beta p(i)+\beta ,\beta }(x_j)\\

&\vdots &\vdots & \vdots &\vdots  &\vdots\\  
&\Phi^{(\beta)}_{\sum\beta p(i),\beta}(x_j) &\Phi^{(\beta )}_{\sum\beta p(i)+\beta_1,\beta}(x_j) & \Phi^{(\beta)}_{\sum\beta p(i)+\beta_2,\beta}(x_j) &\cdots &\Phi^{(\beta )}_{\sum\beta p(i)+\beta,\beta}(x_j)
\end{array}
  \right].
\end{equation}

\subsection{Multivariate Interpolation preserving linear systems of PDEs }
Similar to the univariate case one can adapt the preceding algorithm to the interpolation of multivariate functions, i.e. interpolate $f$ by a polynomial $p$ such that $f=p$, and
\begin{equation}
L_if(x)=L_ip(x) \mbox{ for } x\in \Theta.
\end{equation}
where $\Theta=\{x_0,\cdots ,x_N\}$ is the set of interpolation points, and the partial differential operators are defined by
\begin{equation}
L_if(x)=\sum_{|\alpha|\leq q_i} a^i_{\alpha}(x)\partial^{\alpha}f(x), ĩ=1,\cdots,k.
\end{equation}
The procedure is analogue to that described in Section 2.2. (cf.also \cite{Ka2}).

\section{Approximation of global solutions of linear partial differential equations}

	We refine the algorithm further in order to solve linear partial differential equations globally. In this case the function $u$ to be approximated is not known.  In this section we shall simply describe an algorithm which constructs a polynomial which satifies a linear system of partial differential equations
on an arbitrary set of interpolation points. It is not clear, however, if this polynomial approximation converges to the solution of the system. To ensure that and in order to estimate the rate of convergence we shall need the a priori estimates and regularity results. Note however, that the regularity constraints on the solution maybe low for problems on compact domains as any continuous solution functions $u$ can be approximated by a families of polynomial functions approximating $u$. Therefore, principally, the families of polynomial functions  constructed here may approximate continuous global solutions in viscosity sense. An investigation of this problem will be considered elsewhere in a more general framework where we include some class of nonlinear problems. In order to make the basic ideas transparent we consider first scalar linear problems. We exemplify our algorithm first in the case of dimension $n=1$ and then generalize to the case $n>2$. What we have in mind here are elliptic equations but we need the ellipticity condition only when we wan to prove that the family of polynomials construxted converges to the global solutions. Then we exemplify our method in the case of a typical linear first order system. It is then clear how to generalize to systems of linear equations of any order.

\subsection{The case scalar second order equations of dimension $n=1$}  

We consider the simple boundary value problem

\begin{equation}\label{bound1dim}
L_1u\equiv a(x)\frac{d^2 u}{dx^2}+b(x)\frac{d u}{dx}+c(x)u=f(x) \mbox{ on } (d,e)\subset {\mathbb R},
\end{equation}
with the boundary condition $u(d)=c_d$ and $u(e)=c_e$ (actually an ordinary differential equation). If $a(x)\geq \lambda >0$ for all $x\in {\mathbb R}$, then we have an elliptic operator, but this is not an assumption which we need to construct an univariate polynomial which satisfies the boundary problem on the interpolation points.

We start with the point $d$. We construct a list of polynomial $q_m, m\geq 0$. We define the $q_m$ in substeps. Let $p_0=a_0$. In order that $p_0$ satisfies the boundary condition at $x=d$ we impose
\begin{equation}
p_0=a_0=c_d
\end{equation}
Next we define
\begin{equation}
p_1(x)=a_0+a_1(x-d)
\end{equation}
In order to satisfy the second boundary condition we get
\begin{equation}
p_1(e)=a_0+a_1(e-d)=c_d+a_1(e-d)=c_e \Rightarrow a_1=\frac{c_e-c_d}{e-d}.
\end{equation}
It is clear that $p_1$ preserves the boundary conditions, i.e. $p(d)=u(d)=c_d$ and $p(e)=u(e)=c_e$.
Next let $x_0$ be the first interpolation point (any point in the interval $\left(d,e\right)$. We want to ensure that
\begin{equation}
a(x_0)\frac{d^2 p}{dx^2}(x_0)+b(x_0)\frac{d p}{dx}(x_0)+c(x_0)p(x_0)=f(x_0).
\end{equation}
In order to ensure this, we define a polynomial which is an extension of $p_0$ in three steps. First, define
\begin{equation}
p_2(x)=a_0+a_1(x-d)+a_4(x-x_0)^2(x-d)(x-e) 
\end{equation}
Plugging in and evaluating at $x=x_0$ we get
\begin{equation}
a(x_0)2a_4(x_0-d)(x_0-e)+b(x_0)a_1+c(x_0)(a_0+a_1(x_0-d))=f(x_0)
\end{equation}
Since $a_0,a_1$ are known we get (recall that $x_0\neq d$ and $x_0\neq e$)
\begin{equation}
a_4=\frac{f(x_0)-c(x_0)(a_0+a_1(x_0-d))-b(x_0)a_1}{2a(x_0)(x_0-d)(x_0-e)}.
\end{equation}
Next define
\begin{equation}
p_3(x)=p_2(x)+a_3(x-x_0)(x-x_d)(x-x_e).  
\end{equation}
Plugging in and evaluating at $x=x_0$ we get (assuming that )
\begin{equation}
\begin{array}{ll}
L_1p_3(x_0)=&
L_1p_2(x_0)+a(x_0)a_3(2(x_0-d)\\
\\
&+2(x_0-x_e))+b(x_0)a_3(x_0-d)(x_0-e) =f(x_0).
\end{array}
\end{equation}
Hence, (provided that $x_0\neq d$ and $x_0\neq e$),
\begin{equation}
a_3=\frac{f(x_0)-L_1p_2(x_0)}{a(x_0)(2(x_0-d)+2(x_0-e))+b(x_0)(x_0-d)(x_0-e)}
\end{equation}
Finally, finishing the first inductive step of recursive definition of the polynomial family $(q_m)_{m\in {\mathbb N}}$ 
\begin{equation}
p_4(x)=p_3(x)+a_2(x-d)(x-e).  
\end{equation}
Plugging in and evaluating at $x=x_0$ we get (assuming that )
\begin{equation}
L_1p_4(x_0)=L_1p_3(x_0)+2a(x_0)a_2+b(x_0)((x_0-d)+(x_0-e))=f(x_0).
\end{equation}
Hence, (recall again that $x_0\neq d$ and $x_0\neq e$),
\begin{equation}
a_2=\frac{f(x_0)-L_1p_3(x_0)-b(x_0)((x_0-d)+(x_0-e))}{2a(x_0)((x_0-d)+2(x_0-e))}
\end{equation}
Now we can define
\begin{equation}
q_1(x)=p_4(x) 
\end{equation}
Next assume that the polynomials $q_1, \cdots, q_k$ have been defined. This means that we have computed the polynomial coefficients $a_0,a_1,\cdots, a_{2+3k}$. Then $q_{k+1}$ is defined via
\begin{equation}
q_{k+1}(x)=q_k(x)+(x-d)^3(x-e)^3\Pi_{l=0}^{k}(x-x_l)^3z_k(x),
\end{equation}
where $z_k$ is a polynomial function which will be defined in three substeps.
First, let 
\begin{equation}
q_{k+1,1}(x)=q_k(x)+a_{2+3(k+1)}(x-x_{k+1})^2(x-d)^3(x-e)^3\Pi_{l=0}^{k}(x-x_l)^3
\end{equation}
Plugging in leads to
\begin{equation}
\begin{array}{ll}
L_1q_{k+1,1}(x_{k+1})=&L_1q_k(x_{k+1})+a(x_{k+1})2a_{2+3(k+1)}(x_{k+1}-d)^3\times\\
\\
&(x_{k+1}-e)^3\Pi_{l=0}^{k}(x_{k+1}-x_l)^3=f(x_{k+1}).
\end{array}
\end{equation}
Hence,
\begin{equation}
a_{2+3(k+1)}=\frac{f(x_{k+1})-L_1q_k(x_{k+1})}{a(x_{k+1})2(x_{k+1}-d)^3(x_{k+1}-e)^3\Pi_{l=0}^{k}(x_{k+1}-x_l)^3}
\end{equation}
Next, let 
\begin{equation}
q_{k+1,2}(x)=q_{k+1,1}(x)+a_{2+3k+2}(x-x_{k+1})(x-d)^3(x-e)^3\Pi_{l=0}^{k}(x-x_l)^3
\end{equation}
We define
\begin{equation}
R(x)=(x-d)^3(x-e)^3\Pi_{l=0}^{k}(x-x_l)^3.
\end{equation}
Plugging in leads to
\begin{equation}
\begin{array}{ll}
L_1q_{k+1,2}(x_{k+1})=L_1q_{k+1,1}(x_{k+1})+\\
\\
a(x_{k+1})2a_{2+3k+2}\frac{d^2}{dx^2}R(x_{k+1})+b(x_{k+1})a_{2+3k+2}\frac{d}{dx}R (x_{k+1})=f(x_{k+1}).
\end{array}
\end{equation}
Hence,
\begin{equation}
a_{2+3k+2}=\frac{f(x_{k+1})-L_1q_{k+1,1}(x_{k+1})}{a(x_{k+1})\frac{d^2}{dx^2}R (x_{k+1})+b(x_{k+1})\frac{d}{dx}R (x_{k+1})}
\end{equation}
Finally, let 
\begin{equation}
\begin{array}{ll}
q_{k+1,3}(x)&=q_{k+1,2}(x)+a_{2+3k+1}(x-d)^3(x-e)^3\Pi_{l=0}^{k}(x-x_l)^3\\
\\
&=a_{2+3k+1}R(x)
\end{array}
\end{equation} 
Plugging in leads to
\begin{equation}
\begin{array}{ll}
L_1q_{k+1,3}(x_{k+1})=L_1q_{k+1,2}(x_{k+1})+
a(x_{k+1})a_{2+3k+1}\frac{d^2}{dx^2}R(x_{k+1})\\
\\
+b(x_{k+1})a_{2+3k+1}\frac{d}{dx}R (x_{k+1})+c(x_{k+1})a_{2+3k+1}R(x_{k+1})=f(x_{k+1}).
\end{array}
\end{equation}
Hence,
\begin{equation}
a_{2+3k+2}=\frac{f(x_{k+1})-L_1q_{k+1,2}(x_{k+1})}{a(x_{k+1})\frac{d^2}{dx^2}R(x_{k+1})
+b(x_{k+1})\frac{d}{dx}R (x_{k+1})+c(x_{k+1})R(x_{k+1})}.
\end{equation}
It is clear how to proceed inductively in order to get a family of interpolation polynomials which satisfy the differential equation on an increasing set of interpolation points. Note,however,that we have not used any structural information about the coefficients at this point. This means that the equation may be ill-posed,and convergence cannot be guaranteed.

\subsection{The case of scalar linear partial differential equations}

For a positive integer $k$ consider an equation of form 
\begin{equation}
L_ku\equiv\sum_{|\alpha|\leq k}a_{\alpha}(x)\frac{\partial^{\alpha}u}{\partial x^{\alpha}}=g(x),
\end{equation}
to be solved on the domain $\Omega$ where
\begin{equation}
u\mbox{\Big |}_{\partial \Omega}=f
\end{equation}
What we have in mind is an elliptic equation f order $k$, but ellipticity is not required in order to describe the algorithm which produces a family of multivariate polynomials which satisfy the equation on a set of interpolation points in $\Omega$. Ellipticity becomes important when we want to show that the family of polynomial converges to the solution of the equation (assuming that there is an unique global solution). For simplicity of notation we consider the case $k=2$, i.e. the situation of \eqref{parab}. Assume that $f\in C^k$ and choose a discrete interpolation set $\Theta_b\subset \partial\Omega$. Then we can apply the extended Newton algorithm of Section 3 in order to produce a polynomial $p_b:{\mathbb R}^n\rightarrow {\mathbb R}$ such that
\begin{equation}
\begin{array}{ll}
p_b(x)=f(x) \mbox{ for all } x\in \Theta_b\\
\\
\frac{\partial p_b}{\partial x^{\alpha}}=\frac{\partial p_b}{\partial x^{\alpha}}\mbox{ for all } \alpha \mbox{ with } |\alpha|\leq l \mbox{ and  }x\in \Theta_b
\end{array}   
\end{equation}
We assume that $\Theta_b=\left\lbrace x_{0b},\cdots ,x_{Mb}\right\rbrace $ with $x_{ib}=(x^1_{ib},\cdots ,x^n_{ib})$ and define
\begin{equation}
\Phi_b(x)=\Pi_{i=0b}^{Mb}\Pi_{j=1}^n(x^j-x^j_{i})^{l+1}.
\end{equation}
Next let $\theta_{int}\subset \Omega \setminus \partial \Omega$ be a set of interpolation points in the interior of $\Omega$. Let
\begin{equation}
\Theta_{int}=\left\lbrace x_0,\cdots ,x_N\right\rbrace .
\end{equation}
We enumerate (case $k=2$) the $q:=\frac{(n+1)n}{2}$ diffusion coefficients $a_{\alpha_1},\cdots ,a_{\alpha_q}$ (arbitrary order), where we assume $\alpha_l=(\alpha_{l1},\alpha_{l2})$ and define first $q$ polynomials $p^{\mbox{{\tiny diff}},l}_0(x), l=1,\cdots,q$. Let
\begin{equation}
p^{\mbox{{\tiny diff}},1}_0(x)=p_b(x)+\Phi_b(x)a_{\alpha_1}(x^{\alpha_{11}}-x_0^{\alpha_{11}})(x^{\alpha_{12}}-x_0^{\alpha_{12}}).
\end{equation}
Then we have 
\begin{equation}
L_2p^{\mbox{{\tiny diff}},1}_0(x_0)=L_2p_b(x_0)+\Phi_b(x_0)(1+\delta_{\alpha_{11}\alpha_{12}})a_{\alpha_1}=f(x_0),
\end{equation}
which leads to
\begin{equation}
a_{\alpha_1}=\frac{f(x_0)-L_2p_b(x_0)}{(1+\delta_{\alpha_{11}\alpha_{12}})\Phi_b(x_0)}
\end{equation}
Having defined $p^{\mbox{{\tiny diff}},1}_0(x),\cdots , p^{\mbox{{\tiny diff}},l}_0(x)$ (and therefore computed $a_{\alpha_1},\cdots ,a_{\alpha_l}$) we 
define
\begin{equation}
p^{\mbox{{\tiny diff}},l+1}_0(x)=p^{\mbox{{\tiny diff}},l}_0(x)+\Phi_b(x)a_{\alpha_{l+1}}(x^{\alpha_{(l+1)1}}-x_0^{\alpha_{(l+1)1}})(x^{\alpha_{(l+1)2}}-x_0^{\alpha_{(l+1)2}}),
\end{equation}
and evaluation leads to
\begin{equation}
a_{\alpha_{l+1}}=\frac{f(x_0)-L_2p^{\mbox{{\tiny diff}},l}_0(p_(x_0)}{(1+\delta_{\alpha_{(l+1)1}\alpha_{(l+1)2}})\Phi_b(x_0)}.
\end{equation}
Proceeding inductively we get a $p^{\mbox{{\tiny diff}},q}_0(x)$ which equals together with its derivatives up to order $l$ the function $f$ and such that the diffusion part of the operator applied to $p^{\mbox{{\tiny diff}},q}_0(x)$ equals $g$ at $x_0$. It is now clear how this procedure can be extended such that an extended polynomial $p_0(x)$ equals together with its derivatives up to order $l$ the function $f$ and such that the total operator applied to $p_0(x)$ equals $g$ at $x_0$. As in Section 3 the ansatz for the interpolation polynomial $p_{\Theta}$ which satisfies the linear equation on the set of interpolation points $\Theta=\left\lbrace x_0,\cdots ,x_N \right\rbrace $ then is 
\begin{equation}
p_{\Theta}(x)=\sum_{i=0}^N \Pi_{j=1}^i\Pi_{k=1}^n(x^k-x^k_{j-1})^3p_i(x),
\end{equation}
where $p_i$ for $i\geq 2$ are then constructed as $p_0$ above.

\subsection{The case of  a linear hyperbolic equation}

We consider the hyperbolic equation mentioned above of the form
\begin{equation}
Lu=f \mbox{ in } \Omega, 
\end{equation}
where
\begin{equation}
Lu\equiv \sum_{ij} h_{ij}\frac{\partial u}{\partial x_ii\partial x_j}+\sum_i\frac{\partial}{\partial x_j}+c(x)u
\end{equation}
and $(h_{ij})$ is a symmetric matrix of signature $(n,1)$, if $\mbox{dim}\Omega =n+1$. Note that the operator $L$ can be transformed into the form
\begin{equation}
Lu\equiv \square u + L_1u,
\end{equation}
where $L_1u$ is some first order differential operator on $\Omega$. We assume the initial conditions
\begin{equation}
u=g \mbox{ and } du=\omega ,
\end{equation}
where $g$ and $\omega$ ($1-form$) are initial data.
It is clear that the algorithm described in the preceding section can be used in the present situation. Later we shall see that energy estimates imply convergence of the scheme.

\section{Further refinements: collocation and parallelization}

Numerical experiments show that the coefficients of the recursively computed polynomials have to be computed with increasing accuracy in order to control effects of the truncation error of the coefficients of the polynomials. In the numerical example below, where we computed a polynomial approximation of degree $74$ of the locally analytic function
\begin{equation}
x\rightarrow \frac{1}{1+x}  
\end{equation}
and its derivatives up to order $3$ on the interval $[0,5.4]$ such effects are not observed. However, if we increase the number of derivatives to be approximated up to order $k=10$ and increase the number of interpolation points, effects of truncation errors can be observed for polynomials of degrees larger than 200. The error increases as $|x|$ becomes large and truncation errors increase. This error can be reduced by a more precise representation of the computational approximation of the real numbers involved in the computation. However, as we point out in this section, we can compute $m$ polynomials $p_1^{\Theta_1},\cdots,p_m^{\Theta_m}$ of degree $N_1,N_2\cdots N_m$ parallel which interpolate a given linear system of partial differential equations on some interpolation sets $\Theta_1,\cdots, \Theta_m$ using our basic algorithm, and then compute one polynomial  $p_{\sum\Theta}$ which interpolates the same linear system of partial differential equations on the set
$\sum\Theta =\Theta_1,\cup\cdots, \cup\Theta_m$. It turns out that this can be in such a way that the truncation error of the resulting polynomial $p_{\sum\Theta}$ is much smaller than in case of a direct extension of one polynomial $p_{\Theta_i}$ using the basic algorithm. We call this method the collocation extension of our basic algorithm. We shall assume that the sets of interpolation points are mutually disjunct, i.e.
\begin{equation}
\Theta_i\cap \Theta_j=\oslash \mbox{ iff } i\neq j.
\end{equation}

It is clear that the computation of the polynomials $p_1^{\Theta_1},\cdots,p_m^{\Theta_m}$ can be done parallel and only the step of synthesizing has to be done non-parallel. Next we describe that step in case of two polynomials for simplicit of notation. Extension to $m>2$ polynomials will be clear from that description. So let $\Theta_1,\Theta_2\subset \Omega  \subset {\mathbb R}^n$ be two discrete finite sets of interpolation points of a linear system of partial differential equations $Lu=f$ to be solved on a domain $\Omega$ and such that $\Theta_1\cap\Theta_2=\oslash$. We write down the polynomial in the univariate case because this simplifies the notation, and the multivariate case is quite similar.  
Then we define a regular polynomial interpolation formula on $\Theta_1\cup\Theta_2$ by 
\begin{equation}
\begin{array}{ll}
 \sum_{j=1}^N\Pi_{k\neq j, }\frac{(x-x^{\Theta_1}_k)^{k+1}}{(x_j^{\Theta_1}-x_k^{\Theta_1})^{k+1}}\Pi_{i=1}^{M}
\frac{(x-x^{\Theta_2}_i)^{k+1}}
{(x^{\Theta_1}_j-x_i^{\Theta_2})^{k+1}}
p_{\Theta_1}(x)\\
 \\
 -\sum_j\Pi_{p\in\{1,2\}~l\neq j}
 (x-x_l^{\Theta_p})^{k+1}
a^j_{i,1}(x-x_j^{\Theta_1})^i\\
\\
 +\sum_{j=1}^M\Pi_{k\neq j}
\frac{(x-x^{\Theta_1}_k)^{k+1}}{(x_j^{\Theta_1}
-x_k^{\Theta_1})^{k+1}}\Pi_{i=1}^{N}
\frac{(x-x^{\Theta_2}_i)^{k+1}}{(x^{\Theta_1 }_j-x_i^{\Theta_2})^{k+1}}p_{\Theta_2}(x)\\
\\
-\sum_j\Pi_{p\in\{1,2\}~l\neq j}
(x-x_l^{\Theta_p})^{k+1}
a^j_{i,2}(x-x_j^{\Theta_2})^i\\
\\
=:\sum_j q^{1j}_{\Theta_2,\Theta_1}(x)p_{\Theta_1}(x)+h_1^a(x)\\
\\
+\sum_j q^{2j}_{\Theta_1,\Theta_2}(x)p_{\Theta_2}(x)+h_2^a(x),
\end{array}
\end{equation}
where the constants $a^j_{i,p},p\in\{1,2\}$ are computed recursively as follows: 
For each $j$ we can define $a^j_{0,p}=0$. If $a_{1,p}^j,\cdots, a_{l-1,p}^j$ are determined, then compute $a_{1,p}^j$ via
\begin{equation}
\sum_{1\leq r\leq l}\left(\begin{array}{cc}l\\r\end{array} \right)
D_x^rq^{pj}_{\Theta_1\Theta_2}(x_j)D^{l-r}_x p_{\Theta_p}(x_j)=D^l_x h^a_p(x_j) 
\end{equation}
for each $j$.
Note that this 'synthesis of polynomials' improves the computational power of our method dramatically. In the example below, where we approximate a simple locally analytic function 
\begin{equation}
x\rightarrow \frac{1}{1+x}
\end{equation}
(with convergence radius $1$) and its derivatives up to the third derivative on the interval $[0,5.4]$ with $19$ interpolation points $\Theta_1=\{k0.3|k=0,\cdots 18\}$ we compute a polynomial of degree $74$ in half a minute on a modest laptop machine. If we want to compute a polynomial which gives the same kind of approximation on the interval $[0,5836,8]$ it will take several weeks. However, using parallelization and synthesis, and using the rough estimate that synthesis takes in average the same time as building the $1024$ basis polynomials of degree $74$ on the intervals $[0,5.4]$ and $[k5.7,(k+1)5.7], k=1,\cdots 1023$ we need $10$ steps of parallel synthesis of pairs of polynomials of cost of a less than a minute to get a regular approximation polynomial which is at least of degree $75776$! 
It is clear fromthe preceding remarks how to extend this to the multivariate case (cf. also \cite{Ka2}).

\section{Convergence of polynomial approximations of global solutions of linear elliptic PDE and error estimates by a priori estimates}

Up to now we just considered (regular) polynomial interpolation on given sets of interpolation points.
In this section we consider standard problems in the theory of linear partial differential equations and derive the convergence of our algorithm and error estimates (as the mesh size of the sets of interpolation points converges to zero).
We start with elliptic equations and then consider hyperbolic problems. Similar results can be obtained for initial-value boundary problems for parabolic equations (since analogous error estimates can be obtained). In this case, however, it turns out that (at least for regular data) a WKB-expansion of the fundamental solution has better convergence properties and error estimates can be obtained by Safanov a priori estimates (cf. \cite{KKS} and \cite{Ka}). We shall consider application of our algorithm to this case in the next section. Note that Since to get an error from simple Taylor expansion in genera, because  the interpolated function is unknown. 
\subsection{Convergence for elliptic equations with regular data}
We consider the Dirichlet problem for elliptic equations, i.e. an equation of the form
\begin{equation}\label{elliptsc}
Lu=\sum_{|\alpha|\leq k}a_{\alpha}(x)\frac{\partial u}{\partial x^{\alpha}}=f(x)
\end{equation}
on a domain $\Omega \subseteq {\mathbb R}^n$.  coefficient functions
\begin{equation}
x\rightarrow a_{\alpha}(x),
\end{equation}
and where $u$ is given on the boundary, i.e. 
\begin{equation}
u{\big|}_{\partial \Omega}=g.
\end{equation}

We consider the classical case where $k=2$ and $\Omega$ is bounded. We assume uniform ellipticity, i.e. there exists a constant $K>0$ such that for all $x\in \Omega$
\begin{equation}\label{ell}
\sum_{ij=1}^n a_{ij}(x)\xi_i\xi_j\geq K|\xi|^2.
\end{equation}

In the classical case Schauder boundary estimates are available. We cite them in the context of a standard existence result. for a scalar function $h$ in $\Omega$ we introduce the norms
\begin{equation}
\|h\|_{k}^{bd}=\sum_{j=0}^k\sum_{|\delta|=j}\|D^{\delta}h\|_0
\end{equation}
where
\begin{equation}
\|h\|_0:=\sup_{x\in \Omega}|h(x)|,
\end{equation}
and
\begin{equation}
\|h\|_{k+\alpha}^{bd}=\|h\|_{k}^{bd}+\sum_{j=0}^k\sum_{|\delta|=j}H^{bd}_{\alpha}\left( D^{\delta}h\right),
\end{equation}
where $H_{\alpha}^{bd}(f)$ is the H\"older coefficient of a given function $f$ in $\Omega$.
We assume that the coefficient functions $x\rightarrow a_{ij}(x)$ (diffusion terms), $x\rightarrow b_i(x)$ (drift terms), the potential term ($x\rightarrow c(x)$), and the right side $x\rightarrow f(x)$ are uniformly H\"older continuous (exponent $\alpha$) such that
\begin{equation}\label{coco}
\|a_{ij}\|
^{bd}_{\alpha}
\leq C, 
\|b_{i}\|^{bd}_{\alpha}\leq C, \|c\|^{bd}_{\alpha}\leq C, \|f\|^{bd}_{\alpha}\leq C
\end{equation}
for some generic constant $C$.
\begin{thm}
Assume that conditions \eqref{ell} and \eqref{coco} hold, and assume that $c\leq 0$. Furthermore, assume that $\partial \Omega$ belongs to $C^{2+\alpha}$ and that $g$ belongs to $C_{2\alpha}^bd$. Then the inequalities
\begin{equation}
\begin{array}{ll}
\|u\|_{2+\alpha}^{bd}&\leq C\left(\|g\|+\|u\|_0 +\|f\|_{\alpha}^{bd} \right)\\
\\
&\leq  C\left(\|g\|+\sup_{\partial \Omega}|g|+C\sup_{\Omega}|f| +\|f\|_{\alpha}^{bd} \right)
\end{array}
\end{equation}
hold. Furthermore there exist a unique solution $u\in C_{2+\alpha}^{bd}$ to the Dirichlet problem. 
\end{thm}

The interpolation polynomial $p_{\Theta}$ described in the preceding section is by construction such that
\begin{equation}\label{elliptsc}
L(u-p_{\Theta})=\sum_{|\alpha|\leq k}a_{\alpha}(x)\frac{\partial (u-p_{\theta})}{\partial x^{\alpha}}=\Delta f(x),
\end{equation}
and
\begin{equation}
u-p_{\Theta}{\big|}_{\partial \Omega}=\Delta g.
\end{equation}
It follows that
\begin{thm} Assume the same conditions as in theorem 6.1.. Then
\begin{equation}
\begin{array}{ll}
\|u-p_{\Theta}\|_{2+\alpha}^{bd}
\leq  C\left(\|\Delta g\|+\sup_{\partial \Omega}|\Delta g|+C\sup_{\Omega}|\Delta f| +\|\Delta f\|_{\alpha}^{bd} \right)
\end{array}
\end{equation}
\end{thm}
Note that this implies an $L^2$-error even for the second derivatives of the global solution function, hence essentially an estimate in $H^2(\Omega)$. Even stronger results can be obtained if additional equations for the derivatives of $u$ are considered (cf. \cite{Ka2}).

\subsection{Convergence for a hyperbolic linear partial differential equations equation}

We consider again the hyperbolic equation mentioned above of the form
\begin{equation}\label{hyperb}
Lu=f \mbox{ on } O\subset \Omega, 
\end{equation}
where
\begin{equation}
Lu\equiv \sum_{ij} h_{ij}\frac{\partial u}{\partial x_ii\partial x_j}+\sum_i\frac{\partial}{\partial x_j}+c(x)u
\end{equation}
and $(h_{ij})$ is a symmetric matrix of signature $(n,1)$, if $\mbox{dim}\Omega =n+1$. We assume that some $O\subset \Omega$ is bounded by two spacelike surfaces$\Sigma_i$ and $\Sigma_e$ and swept out by a family of spacelike surfaces $\Sigma_e(s)$. Recall that the initial conditions
\begin{equation}\label{boundary}
u=g \mbox{ and } du=\omega .
\end{equation}

Let $p$ be the interpolation polynom described above such that
\begin{equation}
L(u-p)=\Delta f \mbox{ on } O\subset \Omega.
\end{equation}
\begin{equation}\label{boundary}
u-p=\Delta g \mbox{ and } du=\Delta \omega .
\end{equation}

Then we use the following energy estimate
\begin{prop}
Let $u$ solve the intial value problem \eqref{hyperb}, \eqref{boundary}. Let
\begin{equation}
O(s)=\overline{O}\cap \left\lbrace t\leq s \right\rbrace 
\end{equation}
(swept out by the spacelike surfaces $\Sigma_e(s)$). Then
\begin{equation}
\begin{array}{ll}
\int_{O(s)}|u|^2dV\leq \\
\\
\int_{\Sigma^b_i(s)} |g|^2dS+C(s-s_0)\int_{\Sigma_i}\left( |g|^2+|\omega|^2 \right) dS+C\int_{O(s)}|f|^2dV
\end{array}
\end{equation}
for $s\in \left[s_0,s_1\right]$.
\end{prop}
This implies
\begin{thm} With the same assumptions as in propostion 6.2. we have 
\begin{equation}
\begin{array}{ll}
\int_{O(s)}|u-p|^2dV\leq \\
\\
\int_{\Sigma^b_i(s)} |\Delta g|^2dS+C(s-s_0)\int_{\Sigma_i}\left( |\Delta g|^2+|\omega|^2 \right) dS+C\int_{O(s)}|\Delta f|^2dV
\end{array}
\end{equation}
for $s\in \left[s_0,s_1\right]$.
\end{thm}
Hence the polynomial interpolation scheme described in Section 4 leads to $L^2$-convergence. One can improve this scheme assuming regularity of solutions and considering systems of equations including equations for derivatives of the solution $u$ (cf. \cite{Ka2}).
\section{Applications to parabolic equations (connection to WKB-expansions)} 

We summarize some results concerning WKB-expansions of parabolic equations (cf. \cite{Ka} for details). 
Let us consider the parabolic  diffusion operator
\begin{equation}\label{PPDE} 
\begin{array}{l}
\frac{\partial u}{\partial  t}-Lu\equiv \frac{\partial u}{\partial t}-\frac{1}{2}\sum_{i,j}a_{ij}\frac{\partial^2 u}{\partial x_i\partial x_j}-
\sum_i b_i\frac{\partial u}{\partial x_i},
\end{array}
\end{equation}
where the diffusion coefficients $a_{ij}
$ and the first order coefficients $b_i$ in \eqref{PPDE}
depend on the spatial variable $x$ only. In the following let $\delta t=T-t$, and let
\begin{equation}
(x,y)\rightarrow d(x,y)\ge0,~~(x,y)\rightarrow c_k(x,y),~k\geq 0
\end{equation}
denote some smooth functions on the domain ${\mathbb R}^n\times {\mathbb R}^n$.
Then a set of (simplified) conditions sufficient for pointwise valid WKB-representations 
of the form
\begin{equation}\label{WKBrep}
p(\delta t,x,y)=\frac{1}{\sqrt{2\pi \delta t}^n}\exp\left(-\frac{d^2(x,y)}{2\delta t}+\sum_{k= 0}^{\infty}c_k(x,y)\delta t^k\right), 
\end{equation}
for the solution $(t,x)\rightarrow p(\delta t,x,y)$.
\begin{equation}
\begin{array}{ll}
\frac{\partial u}{\partial  \delta t}- Lu =0, \mbox{with final value}\\
\\
u(0,x,y)=\delta(x-y),
\end{array}
\end{equation}
is given by
\begin{itemize}
\item[(A)] The operator $L$ is uniformly elliptic in ${\mathbb R}^n$, i.e. the matrix norm of $(a_{ij}(x))$ is bounded below and above by $0<\lambda <\Lambda <\infty$ uniformly in~$x$,

\item[(B)] the smooth functions $x\rightarrow a_{ij}(x)$ and $x\rightarrow b_i(x)$ and all their derivatives are bounded.

\end{itemize}
For more subtle (and partially weaker conditions) we refer to \cite{Ka}. We consider the case where there exists a global transformation to the Laplace operator.
If we add the uniform boundedness condition
\begin{itemize}
\item[(C)] there exists a constant $c$ such that for each multiindex $\alpha$ and for all $1\leq i,j,k\leq n,$ 
\begin{equation}\label{unibd}
{\Big |}\frac{\partial
a_{jk}}{\partial x^{\alpha}}{\Big |},~{\Big |}\frac{\partial
b_{i}}{\partial x^{\alpha}}{\Big |}\leq c\exp\left(c|x|^2 \right), 
\end{equation}  
\end{itemize}
then the function $d^2=(x-y)^2$ (in the transformed coordinates and $c_k$ equals its Taylor expansion around $y\in {\mathbb R}^n$, i.e $c_k,k\geq 0$ have the power series representations
\begin{equation}
\begin{array}{ll}
c_k(x,y)&=\sum_{\alpha}c_{k,\alpha}(y)\delta x^{\alpha}, k\geq 0.
\end{array}
\end{equation}
Moreover $c_k, k\geq 0$ are determined by the recursive equations 

\begin{equation}\label{c01e}
-\frac{n}{2}+\frac{1}{2}Ld^2+\frac{1}{2}\sum_{i} \left( \sum_j\left( a_{ij}(x)+a_{ji}(x)\right) \frac{d^2_{x_j}}{2}\right) \frac{\partial c_{0}}{\partial x_i}(x,y)=0,
\end{equation}
where the boundary condition 
\begin{equation}\label{c01b}
c_0(y,y)=-\frac{1}{2}\ln \sqrt{\mbox{det}\left(a^{ij}(y) \right) }
\end{equation}
determines $c_0$ uniquely for each $y\in {\mathbb R}^n$, and
for $k+1\geq 1$ we have
\begin{equation}\label{1gaa}
\begin{array}{ll}
(k+1)c_{k+1}(x,y)+\frac{1}{2}\sum_{ij} a_{ij}(x)\Big(
\frac{d^2_{x_i}}{2}\frac{\partial c_{k+1}}{\partial x_j}
+\frac{d^2_{x_j}}{2} \frac{\partial c_{k+1}}{\partial x_i}\Big)\\
\\
=\frac{1}{2}\sum_{ij}a_{ij}(x)\sum_{l=0}^{k}\frac{\partial c_l}{\partial x_i} \frac{\partial c_{k-l}}{\partial x_j}
+\frac{1}{2}\sum_{ij}a_{ij}(x)\frac{\partial^2 c_k}{\partial x_i\partial x_j}  
+\sum_i b_i(x)\frac{\partial c_{k}}{\partial x_i},
\end{array}
\end{equation}
with boundary conditions
\begin{equation}\label{Rk}
c_{k+1}(x,y)=R_k(y,y) \mbox{ if }~~x=y,
\end{equation}
$R_k$ being the right side of \eqref{1gaa}. In case $a_{ij}=\delta_{ij}$
we have the representations
\begin{equation}\label{eqd}
d^2(x,y)=\sum_i (x_i-y_i)^2,
\end{equation} 

\begin{equation}\label{solc0}
c_0(x,y)=\sum_i(y_i-x_i)\int_0^1 b_i(y+s(x-y))ds,
\end{equation}
and
\begin{equation}\label{solck}
c_{k+1}(x,y)=\int_0^1 R_k(y+s(x-y),y)s^{k}ds,
\end{equation}
$R_k$ being again the right-hand-side of \eqref{1gaa}. The integrals can be taken out if the functions $x\rightarrow b_i(x)$ are given by multivariate power series and error estimates for the truncation error in space and time are obtained (cf. \cite{Ka, KKS}). However, even if the coefficient functions are analytic, i.e. equal locally a power series, it is not possible to approximate such a function globally by their Taylor polynomial. As an example consider the equation
\begin{equation}
\frac{\partial u}{\partial t}-\frac{1}{2}\Delta u -\sum_i^n \frac{1}{1+x_i}\frac{\partial u}{\partial x_i}=0      
\end{equation}
Here, the coefficient functions
\begin{equation}\label{coeff}
x_i\rightarrow \frac{1}{1+x_i} =b_i(x)
\end{equation}
are univariate locally analytic function with convergence radius $1$. Such type of equations occur in praxis of finance (cf. \cite{FrKa, KKS}). In order to obtain an approximation of the WKB-expansion say up to order $5$, i.e. compute the coefficient functions 
\begin{equation}
x\rightarrow c_k(x,y), k=0,\cdots, 5 ,
\end{equation}
we need a global approximation of the functions \eqref{coeff} and their derivatives up to order 10! This is due to the recursion equations for the $c_k, k\geq 1$ which involve second derivatives of $c_{k-1}$. If we have $20$ interpolation points on the $x$-axis this implies that our regular interpolation algorithm computes a polynomial of order $231$. We do the computation in a more modest example in order to keep the resulting polynomial representable on one page in the following section.
\section{A numerical example}
The following polynomial is a similtaneous approximation of the function
\begin{equation}
\begin{array}{ll}
f: [0,5,4]\subseteq {\mathbb R}\rightarrow {\mathbb R}\\
\\
f(x)=\frac{1}{1+x}
\end{array} 
\end{equation}
and its first, second, and third derivative on the domain $[0,5,4]$ with $19$ interpolation points. Hence the degree of this univariate polynomial is $74$. Note that the convergence radius of $f$ is $1$.
\begin{equation}
p_{76}(x)=\sum_{m=0}^{75}a_m(x-x_{m \mbox{{\tiny div}}4})^{m \mbox{{\tiny mod}}4}\Pi_{l=0}^{{ m \mbox{{\tiny div}} 4 -1}}(x-x_l)^{4}
\end{equation}
Note that
\begin{equation}
\frac{d^n}{d x^n}\left(\frac{1}{1+x}\right)|_{x=0}=\frac{(-1)^n n!}{(1+x)^{n+1}}|_{x=0}=(-1)^n n! 
\end{equation}
This leads to the values $a_0=1$, $a_1=-1$, $a_2=1$, and $a_3=-1$ for the coefficients of our interpolation polynomial at $x_0=0$. Note that the coefficients $a_i$ of the interpolation polynomial tend to become smaller for large indexes $i$ as you would expect.


\begin{equation}
\begin{array}{ll}
a_0=1.0\\
\\
a_1=-1.000000000000,\hspace{0.1cm}
a_2=1.000000000000,\hspace{0.1cm}
a_3=-1.000000000000\\
\\
a_4=0.769230765432, \hspace{0.1cm}
a_5=-0.591715921811,\hspace{0.1cm}
a_6=0.455165657066\\
\\
a_7=-0.350124490177,\hspace{0.1cm}
a_8=0.218822618520,\hspace{0.1cm}
a_9=-0.136753442090\\
\\
a_{10}=0.085452696109,\hspace{0.1cm}
a_{11}=-0.053404312398,\hspace{0.1cm}
a_{12}=0.028144617397\\
\\
a_{13}=-0.014953338935,\hspace{0.1cm}
a_{14}=0.008262188243,\hspace{0.1cm}
a_{15}=-0.005370784216\\
\\
a_{16}=0.003734873988,\hspace{0.1cm}
a_{17}=-0.003633027430,\hspace{0.1cm}
a_{18}=0.004645502211\\
\\
a_{19}=-0.006813570816,\hspace{0.1cm}
a_{20}=0.007086610952,\hspace{0.1cm}
a_{21}=-0.007144432312\\
\\
a_{22}=0.006238342564,\hspace{0.1cm}
a_{23}=-0.002646146059,\hspace{0.1cm}
a_{24}=-0.002374282360\\
\\
a_{25}=0.008387675067,\hspace{0.1cm}
a_{26}=-0.015766978592,\hspace{0.1cm}
a_{27}=0.024857498610\\
\\
a_{28}=-0.025373687351,\hspace{0.1cm}
a_{29}=0.025025735340,\hspace{0.1cm}
a_{30}=-0.023974098174\\
\\
a_{31}=0.022321168853,\hspace{0.1cm}
a_{32}=-0.015945627926,\hspace{0.1cm}
a_{33}=0.011155207224\\
\\
a_{34}=-0.007619506803,\hspace{0.1cm}
a_{35}=0.005069120726,\hspace{0.1cm}
a_{36}=-0.002759684498\\
\\
a_{37}=0.001479716734,\hspace{0.1cm}
a_{38}=-0.000790172686,\hspace{0.1cm}
a_{39}=0.000430223475\\
\\
a_{40}=-0.000208511304,\hspace{0.1cm}
a_{41}=0.000106279314,\hspace{0.1cm}
a_{42}=-0.000056281013\\
\\
a_{43}=0.000028862889,\hspace{0.1cm}
a_{44}=-0.000011153733,\hspace{0.1cm}
a_{45}=0.000002201139\\
\\
a_{46}=0.000002233629,\hspace{0.1cm}
a_{47}=-0.000004202246,\hspace{0.1cm}
a_{48}=0.000003699371\\
\\
a_{49}=-0.000002870941,\hspace{0.1cm}
a_{50}=0.000002068390,\hspace{0.1cm}
a_{51}=-0.000001402599\\
\\
a_{52}=0.000000753699,\hspace{0.1cm}
a_{53}=-0.000000375935,\hspace{0.1cm}
a_{54}=0.000000159621\\
\\
a_{55}=-0.000000037499,\hspace{0.1cm}
a_{56}=-0.000000015690,\hspace{0.1cm}
a_{57}=0.000000032004\\
\\
a_{58}=-0.000000031616,\hspace{0.1cm}
a_{59}=0.000000023406,\hspace{0.1cm}
a_{60}=-0.000000010968\\
\\
a_{61}=0.000000001564,\hspace{0.1cm}
a_{62}=0.000000005590,\hspace{0.1cm}
a_{63}=-0.000000011521\\
\\
a_{64}=0.000000012190,\hspace{0.1cm}
a_{65}=-0.000000012095,\hspace{0.1cm}
a_{66}=0.000000011769\\
\\
a_{67}=-0.000000011467,\hspace{0.1cm}
a_{68}=0.000000008698,\hspace{0.1cm}
a_{69}=-0.000000006652\\
\\
a_{70}=0.000000005132,\hspace{0.1cm}
a_{71}=-0.000000003988,\hspace{0.1cm}
a_{72}=0.000000002523\\
\\
a_{73}=-0.000000001600,\hspace{0.1cm}
a_{74}=0.000000001015,\hspace{0.1cm}
a_{75}=-0.000000000643
\end{array}
\end{equation}

\section{Conclusion}

We have designed regular polynomial interpolation algorithms and variations which produce families of multivariate polynomials which solve linear systems of partial differential equations on arbitrary sets of interpolation points. In our basic algorithm the members of the family of polynomials are defined recursively each being an extension of the preceding member in the sense that the preceding member agrees with a given member on the set of interpolation points on which the preceding member satisfies the linear system of partial differential equations. We have shown that the family of multivariate polynomials has the global solution as its natural limit if some a priori information on the system of partial differential equations is available. The information needed can variate from case to case. In any case a solution should exist. We have shown how to use a priori estimates of elliptic equations and of hyperbolic systems of equations in order to obtain error estimates adapted to the regularity of the solution. Similar is true for parabolic equations. All this makes our approach compatible with new techniques like sparse grids or weighted Monte-Carlo algorithms developed in order to treat systems of higher dimension. In case of parabolic equations we showed how regular polynomial interpolation of known functions can be used in order to compute higher order approximations of WKB-expansions of fundamental solutions. We also constructed extensions where the algorithm is parallelized on different set of interpolation points an showed how these partial polynomial approximations can be patched together to one multivariate polynom which fits the given system of linear partial differential equations on the union of sets of interpolation points.

\end{document}